\documentclass[12pt]{article}
\usepackage{fullpage,graphicx,psfrag,amsmath,amsfonts,verbatim}
\usepackage[small,bf]{caption}\usepackage{amssymb}
\usepackage{fullpage}
\usepackage{enumitem}
\usepackage{tikz}
\usepackage{listings}
\usepackage{authblk}

\graphicspath{{figures/}}

\DeclareFixedFont{\ttb}{T1}{txtt}{bx}{n}{12} 
\DeclareFixedFont{\ttm}{T1}{txtt}{m}{n}{12}  

\usepackage{color}
\definecolor{deepblue}{rgb}{0,0,0.5}
\definecolor{deepgreen}{rgb}{0,0.5,0}

\usepackage{hyperref}
\hypersetup{
colorlinks   = true,    
urlcolor     = blue,    
linkcolor    = blue,    
citecolor    = red      
}

\newcommand{\ones}{\mathbf 1}
\newcommand{\reals}{{\mbox{\bf R}}}

\newcommand{\diag}{\mathop{\bf diag}}

\newcommand{\Expect}{\mathop{\bf E{}}}
\newcommand{\Prob}{\mathop{\bf Prob}}

\newcommand{\epi}{\mathop{\bf epi}} 

\newcommand{\eg}{{\it e.g.}}
\newcommand{\ie}{{\it i.e.}}

\newcommand{\BEAS}{\begin{eqnarray*}}
\newcommand{\EEAS}{\end{eqnarray*}}
\newcommand{\BEA}{\begin{eqnarray}}
\newcommand{\EEA}{\end{eqnarray}}
\newcommand{\BEQ}{\begin{equation}}
\newcommand{\EEQ}{\end{equation}}
\newcommand{\BIT}{\begin{itemize}}
\newcommand{\EIT}{\end{itemize}}

\newcommand{\todo}[1]{\textcolor{blue}{\textbf{TODO: #1}}}

\title{Portfolio Construction with Gaussian 
Mixture Returns and Exponential Utility via Convex Optimization}
\author[1]{Eric Luxenberg}
\author[1]{Stephen Boyd}
\affil[1]{Department of Electrical Engineering, Stanford University}

\begin{document}
\maketitle

\begin{abstract}
We consider the problem of choosing an optimal portfolio, assuming the asset 
returns have a Gaussian mixture (GM) distribution, with the objective of 
maximizing expected exponential utility.
In this paper we show that this problem is convex, and readily solved exactly using 
domain-specific languages
for convex optimization, without the need for sampling or scenarios.
We then show how the closely related problem of minimizing 
entropic value at risk can also be formulated as a convex optimization
problem.
\end{abstract}

\newpage
\tableofcontents
\newpage

\section{Introduction}

\subsection{Asset return distributions}
There is a long history of researchers observing that
the tails of asset returns are not well modeled by a Gaussian distribution,
going back to the thesis of Fama in 1965 \cite{Fama},
who observed that
while somewhat symmetric, the tails of the return distribution were much
heavier than those of a Gaussian distribution. Additionally, asset returns are
skewed, violating normality \cite{neuberger-skewness}. 
There is also long history of researchers proposing alternative distributions
to model asset returns, including elliptical distributions \cite{bingham-2001},
and Gaussian mixtures (GMs) \cite{ball-1983,akgiray-1987}, the focus of this paper.

GMs can in principle approximate any continuous 
distribution; for asset returns, it has been observed that good
approximations can be obtained with just a handful of mixture 
components \cite{kon-1984}.
We can interpret the components of a GM return distribution as market regimes, 
with a latent variable that represents the active regime, and
a return distribution that is Gaussian, given the regime.
Many authors have observed that the
correlations among asset returns can change during different market regimes,
for example, increased correlations during bear markets
\cite{bear-correlation,ang2004regimes,ang2015regimes}.  A GM can model such regime-dependent
correlation structures. Another desirable atribute of a GM is that it can model
skewness in return distributions, for which many authors have argued real investors
exhibit preferences \cite{arditti1967risk,scott1980direction}.

GM return models arise in a hidden Markov Gaussian model \cite{ryden-1998} of returns,
which models the regimes as Markovian, and the returns as 
Gaussian, given the regime. 1/
In such a model the means and covariances
of the Gaussian components corresponding to the regimes are fixed, but the 
component weights change in each period \cite{gupta-2012}; but in each period, the 
asset return distribution, conditioned on the past returns, is GM, so the methods
in this paper can be applied.

\subsection{Mean-variance versus expected utility}

\paragraph{Mean-variance portfolio construction.}
In mean-variance portfolio construction, pioneered by Markowitz \cite{Markowitz1952},
portfolio constuction is viewed as an optimization problem with 
two main objectives: the mean or expected return of the portfolio, and the risk, taken
to be the variance of the portfolio return.  These objectives are
combined into a risk-adjusted return using a positive weight parameter, 
interpreted as setting the level of risk aversion.
Mean-variance portfolio construction can be carried out
analytically, when there are very simple constraints,
or numerically, with realistic constraints,
by solving a convex optimization problem such as a 
quadratic program \cite{grinold-1999,cvxbook,boyd-2017,osqp}.
With current convex optimization methods, mean-variance construction 
can be done reliably 
and quickly for up to thousands of assets, and many more
when a factor risk model is used.
These optimization problems can be solved in well under one second, 
allowing back-tests (what-if simulations, based on real or 
simulated data) to be carried out quickly \cite{boyd-2017,cvxpygen,mathworks}.

One obvious criticism of mean-variance portfolio construction 
is that the (quadratic) objective function penalizes returns that are well
above the mean (a desirable event) just as much as returns
that are well below the mean (an undesirable event) \cite{hanoch1970efficient}.
Another is that it only uses the first two moments of the 
return distribution, and so cannot take into account skewed or 
fat-tailed distributions.  Nevertheless it is very widely used
in practice. There is work in analyzing portfolios under higher order Taylor
approximations of utilities \cite{jondeau2006optimal}, but these are not used in practice due to
both semantic and computational complexity. 

\paragraph{Expected utility portfolio construction.}
In work that predates mean-variance portfolio construction,
Von Neumann and Morgenstern \cite{utility} introduced the notion of utility to
model decision making with uncertain outcomes. A utility function specifies a
value indexing an agent's preference for each specific outcome; their
theory posits that the agent makes a choice so as to maximize her expected utility.

Portfolio construction by expected utility maximization also frames the 
problem as an optimization problem.  The trader specifies a utility 
function that is concave and increasing, and the objective (to be 
maximized) is the expected utility under the return distribution.
This formulation avoids the awkward situation in mean-variance portfolio
construction where high portfolio returns are 
considered bad.  Expected utility maxmization better captures the asymmetry in downside 
versus upside risks than mean-variance optimization.  Since
the return distribution is arbitrary, expected utility can directly handle 
return distributions with skew or fat tails.

Expected utility maximization, like mean-variance optimization,
leads to a convex optimization problem,
more specifically, a stochastic optimization problem \cite{shapiro}.
Almost all expected utility methods for constructing portfolios 
work with samples of the asset returns. 
This can be considered an advantage, since it means that such methods
can work with any return distribution from which we can sample returns.
The disadvantage is that sample-based optimization, while tractable,
can be slow compared to mean-variance methods, and scales poorly with problem size.

There are several related portfolio construction methods that rely on
return samples and stochastic convex optimization.
One is based on conditional value at risk (CVaR) \cite{Rockafellar-CVaR,CVaR-gen}. 
A more recently proposed method uses entropic value at risk 
(EVaR) \cite{EVaR-theory,Cajas}, which we address in \S\ref{s-EVaR}.
Both of these are coherent measures of risk \cite{Artzner,rockafellar-coherent},
and result in convex stochastic optimization problems.

Sample based stochastic optimization methods are used in practice \cite{grinold1999mean},
but far less often than methods based on mean-variance optimization,
which do not involve samples.  This is part because solving sample based stochastic 
convex optimization problems is tractable, but far more involved than solving
the convex optimization problems that do not involve return samples, 
\eg, mean-variance optimization or the methods proposed in this paper.

\paragraph{Comparison.}
These two main approaches, mean-variance optimization and expected utility maximization,
are not as different as they might seem. 
Levy and Markowitz \cite{levy-mark-utility} show that maximizing a second order
Taylor approximation of a utility function is equivalent to mean-variance 
optimization.  So very roughly speaking, mean-variance optimization is the 
second order approximation of expected utility optimization.

When the returns are Gaussian, and we use an exponential utility,
mean-variance and expected utility optimization are
not merely close, but exactly the same.
It appears that Merton was the first to note this connection
\cite{Merton1969}, but his observation does not seem to be
mentioned often after that (see also \S\ref{s-special}).

\subsection{This paper}
In this paper we consider a GM model for asset returns,
and maximize expected utility with a generic utility function,
the exponential utility \cite{saha-1993}.
We refer to this type of portfolio construction as EGM,
for exponential utility with Gaussian mixture returns.
We show that the EGM portfolio construction problem can be solved
exactly as a convex optimization problem, without the need for 
any samples from the distribution or other approximations.
The EGM portfolio construction problem is not only convex, 
but is easily specified in just a few lines
of code in domain-specific languages (DSLs) for convex optimization
such as CVXPY \cite{diamond-2016}, CVX \cite{grant2014cvx}, 
or CVXR \cite{cvxr2020}.
Thus EGM combines the efficiency, reliability, and scalability
of mean-variance optimization with the ability of expected utility 
maximization to handle non-Gaussian returns and the asymmetry
in our preferences.
When the GM has only one component, our return model is Gaussian,
and EGM reduces to mean-variance optimization.  Thus we can think of EGM as an 
extension of mean-variance optimization, or as a special case 
of expected utility where the 
problem can be solved exactly, without any return samples.

We also show that EGM is closely related to portfolio construction methods
based on the entropic value at risk (EVaR).  With GM return model,
we show that EVaR portfolio construction problem leads to a convex 
optimization problem that, like EGM, does not involve sample based 
stochastic optimization.

\subsection{Previous and related work}
\paragraph{Portfolio construction with Gaussian mixture returns.} 
In \cite{buckley-gmm},
Buckley et al.\ consider a two-component Gaussian mixture
of tranquil and distressed regimes, and analyze several objectives, including
Markowitz, Sharpe ratio, exponential utility, and lower partial moments.
In \S3.3.6 of their paper
they derive the closed form expression for expected
exponential utility under Gaussian mixture returns,
but do not observe that maximizing exponential utility leads to a convex problem. 
Studying single period portfolios consisting of a risk free asset and a risky
asset, Prigent and Kaffel analyze optimal portfolios under arbitrary utility
functions, and show on historical data that GM return models lead to
significantly different portfolios than those arising from a Gaussian
return model \cite{kaffel-2014}.
\paragraph{EVaR portfolio construction.} In recent work \cite{Cajas},
Cajas develops a disciplined convex
(DCP) formulation of EVaR, with return samples,
which allows it to be used as
either the objective or as a constraint in
portfolio optimization problems specified using DSLs for convex optimization
such as CVXPY.
Since a return distribution that takes on a finite set of values
(\eg, the empirical distribution of samples) is a special case
of GM, we can consider EGM (with EVaR) as a generalization of 
Cajas' formulation.

\subsection{Outline}
We describe the GM return model is \S\ref{s-gmm-ret},
and in \S\ref{s-port-opt} we show that portfolio optimization with 
exponential utility is a convex optimization problem.
In \S\ref{s-EVaR} we show that the closely related function EVaR is also convex, 
so minimizing it, or adding a limit on it as a constraint,
results in a convex optimization problem.

\section{Gaussian mixture return model}\label{s-gmm-ret}

\subsection{Asset return distribution}
We let $r\in \reals^n$ denote the return of $n$ assets over some specific period.
We model $r$ as having a GM distribution with $k$ components,
\[
r\sim \textbf{GM}(\{\mu_i,\Sigma_i,\pi_i\}_{i=1}^k),
\]
where $\pi_i\in \reals$ are the (positive) component probabilities, 
$\mu_i \in \reals^n$ are the component means, and
$\Sigma_i \in \reals^{n \times n}$ are the (symmetric positive definite) 
component covariance matrices.

The GM return distribution includes two interesting special cases.  When there is
only one component, it reduces to Gaussian, with $r \sim \mathcal N(\mu_1, \Sigma_1)$.
Another special case arises when $\Sigma_1 = \cdots = \Sigma_k=0$.  Here
$r$ takes on only the values $\mu_1, \ldots, \mu_k$, with probabilities
$\pi_1, \ldots, \pi_k$. We refer to this as a finite values return distribution.

\subsection{Portfolio return distribution}
Let $w\in\reals^{n}$ denote the weights in an investment portfolio, with 
$\ones^Tw =1$, where $\ones$ is the vector with all entries one. For $w_i\geq 0$,
$w_i$ is the fraction of the total portfolio value invested in asset $i$; for $w_i<0$,
$-w_i$ is the fraction of total portfolio value that is held in a short position in
asset $i$.
The portfolio return is $R=w^Tr$.  This scalar random variable is also GM with
component probabilities $\pi_i$, and means and variances
\[
\nu_i = w^T \mu_i, \quad \sigma_i^2 = w^T\Sigma_i w,  \quad i=1, \ldots, k.
\]

We observe that various quantities associated with
the portfolio return $R$ can be evaluated analytically, 
without the need for Monte Carlo or other sampling methods.    
For example its cumulative distribution
function (CDF) is given by
\BEQ\label{e-gmm-cdf}
\Phi_R(w,a) = \sum_{i=1}^k \pi_i \Phi \left(\frac{a-\nu_i}{\sigma_i}\right) =
\sum_{i=1}^k \pi_i \Phi \left(\frac{a-w^T\mu_i}{(w^T \Sigma_i w)^{1/2}}\right),
\EEQ
where $\Phi$ is the CDF of a standard Gaussian.

\subsection{Moment and cumulant generating functions}
Two other quantities we will encounter later are the moment generating function
\BEQ\label{e-mgf}
M(w,t) = \Expect\exp (t R) =
\sum_{i=1}^k\pi_i \exp \left( t \nu_i+\frac{t^2}{2}\sigma_i^2\right) =
\sum_{i=1}^k\pi_i \exp \left( t \mu_i^Tw+\frac{t^2}{2}w^T\Sigma_i w \right),
\EEQ
where we use $\Expect \exp Z = \exp (\mu+\sigma^2/2)$ for $Z \sim \mathcal N(\mu,\sigma^2)$,
and the cumulant generating function
\begin{eqnarray}
\nonumber
K(w,t) &=& \log \Expect\exp (tR)\\
\nonumber &=& \log \left(
\sum_{i=1}^k\pi_i \exp \left( t\nu_i+\frac{t^2}{2}\sigma_i^2\right)\right)\\
&=&
\log \left( 
\sum_{i=1}^k\pi_i \exp \left( t\mu_i^Tw+\frac{t^2}{2}w^T\Sigma_i w \right)
\right).
\label{e-cgf}
\end{eqnarray}
We observe for future use the identity 
\BEQ\label{e-ident}
K(w,t)=K(tw,1),
\EEQ
\ie, the parameter $t$ simply multiplies the argument $w$.

\section{Portfolio optimization with exponential utility}\label{s-port-opt}

\subsection{Expected exponential utility}
Our objective is to choose $w$ to maximize the expected exponential utility 
$\Expect U_\gamma(R)$, where 
\[
U_\gamma(a) = 1-\exp(-\gamma a),
\]
with $\gamma>0$ the risk aversion parameter.
Using \eqref{e-mgf}, we can express this as 
\[
\Expect U_\gamma(R) = 1-\Expect \exp(-\gamma w^Tr) 
= 1-  M(w,-\gamma).
\]
It follows that we can maximize $\Expect U_\gamma(R)$ by minimizing
the moment generating function $M(w,-\gamma)$,
or equivalently the cumulant generating function
\BEQ\label{e-obj}
K(w,-\gamma) = \log \left( \sum_{i=1}^k\exp \left( \log\pi_i -\gamma\mu_i^Tw+
\frac{\gamma^2}{2}w^T\Sigma_i w \right) \right).
\EEQ

\paragraph{Convexity.}
The function $K(w,-\gamma)$ is a convex function of $w$.  
To see this, we note that for each $i$,
$\log \pi_i- \gamma \mu_i^Tw + \frac{\gamma^2}{2} w^T \Sigma_i w$ is a convex
quadratic function of $w$, and therefore convex.
The function $K(w,-\gamma)$ is the log-sum-exp function (also called the soft-max
function),
\begin{equation}\label{e-log-sum-exp}
S(u) = \log\left(\sum_{i=1}^k\exp u_i\right),
\end{equation}
of these arguments. 
The log-sum-exp function is convex and increasing in all arguments, so 
the composition $f_\gamma$ is convex \cite[\S3.1.5]{cvxbook}.

\subsection{EGM portfolio construction}
Our portfolio construction optimization problem has the form
\begin{equation}\label{e-prob}
\begin{array}{ll}
\mbox{minimize} &  K(w,-\gamma)\\
\mbox{subject to} & \ones^Tw =1, \quad w \in \mathcal W,
\end{array}
\end{equation}
where $\mathcal W$ is a convex set of portfolio constraints.
This is evidently a convex optimization problem.
One implication is that we can efficiently solve this problem globally using
a variety of methods.

\paragraph{DSL specification.}
The EGM problem \eqref{e-prob} is not just convex.
It is readily specified in domain-specific languages (DSLs) for convex
optimization, since all such systems include the log-sum-exp function,
and all such systems can handle the convex function composition rules
that we used to establish convexity of $K(w,-\gamma)$ in $w$.
No special methods (or gradient or other derivatives) are needed; 
the function $K(w,-\gamma)$ can be specified in a DSL by just typing it 
in as is.
As a simple example, CVXPY code for specifying
the EGM construction problem \eqref{e-prob} and solving it,
with a long only portfolio (\ie, $w \geq 0$), is given below.
(This code snippet is also available at the repository 
\url{https://github.com/cvxgrp/exp_util_gm_portfolio_opt}.)

\bigskip

\begin{lstlisting}
import cvxpy as cvx

def K(w):
  u = cvx.vstack([cvx.log(pi[i])
    - gamma * mus[i] @ w
    + (gamma**2/2) * cvx.quad_form(w, Sigmas[i]) for i in range(len(pi))])
  return cvx.log_sum_exp(u)

w = cvx.Variable(n)
objective = cvx.Minimize(K(w))
constraints = [ w >= 0, cvx.sum(w) == 1 ]
egm_prob = cvx.Problem(objective, constraints)
egm_prob.solve()
w.value
\end{lstlisting}

\bigskip

Here it is assumed that \verb|n|, \verb|pi|, and \verb|gamma| are 
constants corresponding to $n$, $\pi$ and $\gamma$,
and \verb|mus| and \verb|Sigmas| are lists 
of the $\mu_i$ and $\Sigma_i$, respectively.
In lines 3--7 the objective $K(w,-\gamma)$ is formed,
and in lines 9--12 the EGM optimization problem is formed.
The problem is solved in line 13, which populates \verb|w.value| with optimal
weights.
In this simple example our portfolio constraint set $\mathcal W$ is simple.
One of the advantages of using a DSL is that more complex constraints 
can be added by just appending them to the list of constraints defined in line~11.

\paragraph{Soft-max interpretation.}
We can give an interpretation of the objective $K(w,-\gamma)$ in \eqref{e-obj} 
in terms of the soft-max function,
which can be thought of as a smooth approximation to the max, since it
satisfies
\BEQ\label{e-soft-max}
\max_i u_i  \leq S(u) \leq \max_i u_i  + \log k.
\EEQ

The objective \eqref{e-obj} can be expressed as
\[
K(w,-\gamma) = S(u), \qquad
u_i = \log \pi_i + \gamma \left( -\mu_i^Tw+ \frac{\gamma}{2}w^T\Sigma_i w\right),
\quad i=1, \ldots, k.
\]
We recognize $-\mu_i^Tw + \frac{\gamma}{2} w^T \Sigma_i w$ as the negative
risk-adjusted return of the portfolio under the $i$th Gaussian component.
Thus $K(w,-\gamma)$ is the soft-max of these negative risk-adjusted returns,
offset by the terms $\log \pi_i$, and scaled by $\gamma$.
Roughly speaking, our objective is an approximation of the maximum of 
the negative risk-adjusted returns under the component distributions.

From \eqref{e-soft-max} we have
\begin{eqnarray}
\label{e-K-lower-bnd}
K(w,-\gamma) &\geq& 
\max_{i=1,\ldots, k} \left( \log \pi_i-\gamma \mu_i^Tw+ \frac{\gamma^2}{2}w^T\Sigma_i w \right),\\
\label{e-K-upper-bnd}
K(w,-\gamma) &\leq& \log k + 
\max_{i=1,\ldots,k}
\left( \log \pi_i-\gamma \mu_i^Tw+ \frac{\gamma^2}{2}w^T\Sigma_i w \right).
\end{eqnarray}

\subsection{Special cases}\label{s-special}
\paragraph{Gaussian returns.}
When $k=1$ our GM return distribution reduces to Gaussian,
and the problem \eqref{e-prob} reduces to the standard Markowitz problem
\cite{Markowitz1952,mark-book} 
\[
\begin{array}{ll}
\mbox{maximize} &  \mu_1^Tw - \frac{\gamma }{2} w^T\Sigma_1 w\\
\mbox{subject to} & \ones^Tw =1, \quad w \in \mathcal W.
\end{array}
\]

\paragraph{Finite values returns.}
When $\Sigma_i=0$, so $r$ takes on only the values $\mu_1, \ldots, \mu_k$,
the problem \eqref{e-prob} can be expressed as
\begin{equation}\label{e-prob-finite-vals}
\begin{array}{ll}
\mbox{minimize} & \log\left( \sum_{i=1}^k \pi_i \exp (- \gamma \mu_i^Tw) \right) \\
\mbox{subject to} & \ones^Tw =1, \quad w \in \mathcal W.
\end{array}
\end{equation}

\subsection{Simple example}
To illustrate the difference between EGM and mean-variance portfolios, 
we construct a very simple example for which
the two portfolios can be analytically found.
We take a finite values distribution with $n=2$ assets and $k=2$ components,
with $\Sigma_1 = \Sigma_2=0$,
\[
\mu_1 = \left[\begin{array}{r} -1 \\ 0 \end{array}\right], \qquad
\mu_2 = \left[\begin{array}{r} 1 \\ 0 \end{array}\right].
\]
Thus $r=(-1,0)$ with probability $\pi_1$ and $r=(1,0)$ with probability $\pi_2$.
The first asset is risky, and the second is riskless, with zero return.
We take $\mathcal W= \reals^2$, so the only constraint on 
the portfolio weight is $w_1+w_2=1$.

\paragraph{Markowitz portfolio.}
The mean and covariance of $r$ are
\[
\mu = \left[\begin{array}{c} 1-2\pi_1 \\ 0 \end{array}\right], \qquad
\Sigma = \left[ \begin{array}{cc} 4\pi_1(1-\pi_1) & 0 \\ 0 & 0 \end{array}\right].
\]
The Markowitz optimal portfolio is 
\[
w^\text{M}_1 = \frac{1-2\pi_1}{4 \gamma \pi_1(1-\pi_1)},
\]
with $w^\text{M}_2 = 1-w^\text{M}_1$.

\paragraph{EGM portfolio.}
The EGM portfolio minimizes
\[
\pi_1 \exp (\gamma w_1) + (1-\pi_1)\exp (-\gamma w_1),
\]
so the EGM portfolio is
\[
w^\text{E}_1 = \frac{\log(1/\pi_1 -1)}{2\gamma},
\]
with $w^\text{E}_2 = 1-w^\text{E}_1$.

\paragraph{Comparison.}
The two portfolios are the same for $\pi_1=1/2$, 
with $w^\text{M} = w^\text{E} =(0,1)$.
They are not too far from each other for other values of $\pi_1$ and $\gamma$,
but can differ substantially for others.
For example with $\pi_1 =0.05$ and $\gamma=1$, 
the Markowitz and EGM portfolios are
\[
w^\text{M}=(4.74,-3.74), \qquad
w^\text{E}=(1.47,-.47).
\]
The value at risk 5\% is 4.74 for the Markowitz porfolio compared
to 1.47 for the EGM portfolio. 

\subsection{High and low risk aversion limits}
\paragraph{High risk aversion limit.}
Here we consider the case where $\gamma \to \infty$.
Dividing \eqref{e-K-lower-bnd} and
\eqref{e-K-upper-bnd} by $\gamma$, we find that
\[
\frac{K(w, -\gamma)}{\gamma} = 
\max_{i=1, \ldots, k} \left(-\mu_i^T w + \frac{\gamma}{2} w^T
\Sigma_i w \right) + O(1/\gamma).
\]
So for large risk aversion parameter $\gamma$, 
the EGM portfolio construction problem \eqref{e-prob} is approximately 
\[
\begin{array}{ll}
\mbox{minimize} &  \max_{i=1, \ldots, k} \left( -\mu_i^T w + \frac{\gamma}{2}
w^T \Sigma_i w \right)\\
\mbox{subject to} & \ones^Tw =1, \quad w \in \mathcal W.
\end{array}
\]
Thus in the limit of high risk aversion,
the EGM portfolio minimizes the maximum of the risk adjusted returns
under each of the components, regardless of the $\pi_i$.
This is similar to solving a minimax Markowitz problem, where we 
use as a risk model the maximum risk over a set of covariance matrices
(see \cite[\S4.2, p.\;30]{boyd-2017}).

\paragraph{Low risk aversion limit.}
Here we consider the case where $\gamma \to 0$.
We start with the well known expansion
\[
\frac{1}{\gamma} \log \left(\sum_{i=1}^k \pi_i \exp \gamma z_i \right)
= \sum_{i=1}^k \pi_i z_i + \frac{\gamma}{2} \left(
\sum_{i=1}^k \pi_i z_i^2 - \left( \sum_{i=1}^k \pi_i z_i\right)^2 \right)
+ O(\gamma^2),
\]
for any $z_i$. (We recognize the first term on the righthand side 
as the mean of $z$, and the 
second as $\gamma/2$ times the variance of $z$, when $z$ is a random variable
taking values $z_1, \ldots, z_k$ with probabilities $\pi_1, \ldots, \pi_k$.)
Substituting $z_i = -\mu_i^T w + \frac{\gamma}{2} w^T \Sigma_i w$ we obtain
\[
\frac{K(w,-\gamma)}{\gamma} = -\mu^T w + \frac{\gamma}{2} w^T \Sigma w + O(\gamma^2),
\]
where $\mu$ and $\Sigma$ are the mean and covariance of $r$, 
\[
\mu = \Expect r = \sum_{i=1}^k \pi_i \mu_i, \qquad
\Sigma = \Expect rr^T -(\Expect r)(\Expect r)^T = 
\sum_{i=1}^k \pi_i \left( \Sigma_i + (\mu_i - \mu)(\mu_i-\mu)^T \right).
\]
This has a very nice interpretation: in the limit of small risk aversion,
EGM reduces to Markowitz, using the mean and covariance of the return.

\section{Portfolio optimization with entropic value at risk} \label{s-EVaR}
\subsection{Entropic value at risk}
The traditional measure of downside risk
is the value at risk (VaR) with probability $\alpha$, which is the $(1-\alpha)$
quantile of the negative return $-R$,
\[
\text{VaR}_\alpha(R) = -\inf\{x\in\reals\mid \Prob(R\leq x)>\alpha\}.
\]
(We are typically interested in values such as $\alpha=0.05$ or $\alpha = 0.01$.)
For example if the value at risk of a portfolio with probability $5\%$ is $15\%$, 
the probability of a loss exceeding $15\%$ (\ie, $R\leq -0.15$) is $5\%$.
Value at risk is interpretable and widely used, 
but it is not a coherent risk measure \cite{Rockafellar-CVaR}.
For example, VaR is not sub-additive, so the sum of two 
portfolios can have a higher VaR than the sum of the component VaRs.
Several coherent risk measures have been proposed, including the
conditional value at risk $\text{CVaR}_\alpha$ \cite{Rockafellar-CVaR}
and entropic value at risk $\text{EVaR}_\alpha$ \cite{EVaR-theory}.

The entropic value at risk $\text{EVaR}_\alpha$ is the tightest 
Chernoff upper bound on $\text{VaR}_\alpha$, which can be expressed in
terms of the cumulant generating function as
\[
\text{EVaR}_{\alpha}(R) = \inf_{\lambda>0} \frac{K(w,-\lambda) -\log \alpha}{\lambda}
\geq \text{VaR}_\alpha(R)
\]
(It is also an upper bound on $\text{CVaR}_\alpha$.) 
Cajas \cite{Cajas} and Shen et al.~\cite{shen2022minimizing}
describe convex optimization problems involving EVaR with the expectation replaced
with its sample approximation. 

\paragraph{Minimum EVaR portfolio.}
To minimize $\text{EVaR}_\alpha(R)$, we solve the 
optimization problem 
\begin{equation}\label{e-EVaR-prob}
\begin{array}{ll}
\mbox{minimize} &  \frac{K(w,-\lambda) -\log \alpha}{\lambda} \\
\mbox{subject to} & \ones^Tw =1, \quad w \in \mathcal W, \quad \lambda>0,
\end{array}
\end{equation}
with variables $w\in \reals^n$ and $\lambda \in \reals$.

There is a close connection between this problem and the 
exponential utility maximization problem \eqref{e-prob}.
Suppose that $w^\star$ and $\lambda^\star$ are optimal for
\eqref{e-EVaR-prob}.  Then $w^\star$ is also optimal for the 
exponential utility problem \eqref{e-prob}, with risk aversion 
parameter $\lambda^\star$.
Thus we can think of minimizing $\text{EVaR}_\alpha(R)$ as 
simply choosing a value of the risk aversion parameter in EGM.
(This choice of parameter depends on $\alpha$.)
We will refer to portfolio construction using \eqref{e-EVaR-prob} also as
EGM, since any such portfolio is optimal for EGM with some
value of risk aversion, and also, conveniently, entropic and 
exponential both start with E.

\subsection{Convex formulation}
The problem \eqref{e-EVaR-prob} is not a convex optimization problem 
since the objective is not jointly convex in $w$ and $\lambda$.
But a change of variable can give us an equivalent convex problem.
Instead of using the variable $\lambda$, we use
the new variable $\delta = 1/\lambda$, and the problem \eqref{e-EVaR-prob}
becomes
\begin{equation}\label{e-EVaR-prob-cvx}
\begin{array}{ll}
\mbox{minimize} &  \delta K(w/\delta,-1) - \delta \log \alpha \\
\mbox{subject to} & \ones^Tw =1, \quad w \in \mathcal W, \quad \delta>0,
\end{array}
\end{equation}
with variables $w\in \reals^n$ and $\delta \in \reals$.
(We use the identity \eqref{e-ident} above.)
This objective is jointly convex in the variables $w$ and $\delta$, since it 
is the perspective function of $K(w,-1)$, which is convex in $w$
\cite[\S 3.2.6]{cvxbook}, so \eqref{e-EVaR-prob-cvx} is a convex 
optimization problem, which is readily solved.
(The constraint $\delta>0$ is actually not needed, since the perspective
function is defined to be $+\infty$ if $\delta \leq 0$.)

Unfortunately, current DSLs for convex optimization do not automate 
the creation of the perspective of a function, so the 
problem \eqref{e-EVaR-prob-cvx} cannot simply be typed in; we must form
the the perspective function by hand, as outlined below 
in appendix~\ref{a-graph-form}.

There are also simple methods that can be used to solve it, with a modest loss
in efficiency, that are immediately compatible with DSLs.
One method is alternating optimization, where we alternate between fixing 
$\delta$ and optimizing over $w$ (easy with current DSLs),
and fixing $w$ and optimizing over $\delta$ (minimization of 
a scalar convex function, which can be done by many simple methods).
To start we can replace $K$ with the lower bound 
\eqref{e-K-lower-bnd} (or the upper bound \eqref{e-K-upper-bnd}), to 
obtain the approximate problem
\begin{equation}\label{e-EVaR-prob-approx}
\begin{array}{ll}
\mbox{minimize} & \max_i \left( \delta \log (\pi_i/\alpha) - \mu_i^T w + 
\frac{w^T\Sigma_i w}{2 \delta} \right) \\
\mbox{subject to} & \ones^Tw =1, \quad w \in \mathcal W, \quad \delta>0.
\end{array}
\end{equation}
This problem is convex, and also immediately representable in DSLs using
the quadratic-over-linear function for the last term in the objective.
(Here too the constraint $\delta>0$ is redundant, since the 
quadratic-over-linear function is defined as $+\infty$ if the 
denominator is not positive.)

\subsection{Special cases}
\paragraph{Gaussian returns.}
When $k=1$, our GM return distribution is Gaussian and we have
\BEQ\label{e-EVaR-gaussian}
\delta K(w/\delta, -1)-\delta \log \alpha =
-\delta \log \alpha - \mu_1^Tw+ \frac{w^T\Sigma_1w}{2 \delta},
\EEQ
with variables $w$ and $\gamma$.  This objective is readily minimized using DSLs,
using the quadratic-over-linear function for the last term.

The value of $\delta$ that minimizes this, with fixed $w$, is 
\[
\delta = \left( \frac{w^T\Sigma_1 w}{-2 \log \alpha}\right)^{1/2}.
\]
Thus we see that for Gaussian returns, the portfolio that minimizes
$\text{EVaR}_\alpha$ is in fact Markowitz, with the specific 
choice of risk aversion parameter
\[
\gamma = \left( \frac{-2 \log \alpha}{w^T\Sigma_1 w} \right)^{1/2}.
\]
(This depends on $w$, so to find it we must solve the convex problem
with objective \eqref{e-EVaR-gaussian}.)
We see that as $\alpha$ decreases, the associated risk aversion increases, which makes sense.

Substituting the optimal value of $\delta$ into the objective \eqref{e-EVaR-gaussian}, 
we find that the objective is
\[
- \mu_1^T w + (-2 \log \alpha)^{1/2} \left( w^T\Sigma_1 w \right)^{1/2},
\]
plus a constant.
Thus we maximize a risk adjusted return, using the standard
deviation instead of the traditional variance as risk,
and the very specific risk aversion constant $(-2\log \alpha)^{1/2}$.

\paragraph{Finite values returns.}
When $\Sigma_i=0$, so $r$ takes on only the values $\mu_1, \ldots, \mu_k$,
the problem \eqref{e-EVaR-prob-cvx} can be expressed as
\begin{equation}\label{e-EVaR-prob-finite-vals}
\begin{array}{ll}
\mbox{minimize} & \delta \log\left( \sum_{i=1}^k (\pi_i/\alpha)
\exp (-\mu_i^Tw/\delta) \right) \\
\mbox{subject to} & \ones^Tw =1, \quad w \in \mathcal W.
\end{array}
\end{equation}

\section{Conclusions}

In this paper we have shown that two specific portfolio construction problems,
maximizing expected exponential utility and minimizing entropic value at risk,
with a Gaussian mixture return model, can be
formulated as convex optimization problems, and exactly solved 
with no need for return samples or Monte Carlo approximations.
The resulting problems 
are not much harder to solve than a mean-variance problem, but have the advantage
of directly handling return distributions with substantial skews or large tails.

Our focus in this paper is on the formulation of the these portfolio 
construction problems as tractable
convex optimization problems that do not need return samples.
In a future paper we will report on practical portfolio construction
using these methods.

\clearpage
\bibliography{library}
\section{Data availability}
Data sharing not applicable to this article as no datasets were generated or analysed during the current study.
\clearpage 
\appendix
%

\section{Graph form representation of EVaR}\label{a-graph-form}

\subsection{Graph form representation}
In this section we show how to express the objective of \eqref{e-EVaR-prob-cvx}
in graph form \cite{grant2008graph}, which is the basic representation of 
a function in DSLs for convex optimization, that rely on disciplined convex programming 
(DCP) \cite{grant2006disciplined}.
In a recent paper Cajas gave a graph form description of EVaR, for the special
case when $\Sigma_i=0$, \ie, for a finite values return model \cite{Cajas}.
Thus we are extending his formulation from a finite values return model 
to a GM return model.

The graph form of a function $f:\reals^n \to \reals$ expresses the 
epigraph of $f$ as the inverse image of a cone under an affine mapping.
(For practical use, the cone must be a Cartesian products of cones supported 
by the solver.)
The graph form of $f$ is
\begin{equation}\label{e-general-epigraph}
\epi f = \{(x,t) \mid f(x) \leq t \} = 
\{ (x,t) \mid\exists z~ Fx+Gz+td+e \in C\},
\end{equation}
where $F\in \reals^{p \times n}$, $G\in \reals^{p \times m}$,
$d\in \reals^p$, and $e\in \reals^p$
are the coefficients, and $C\subseteq \reals^p$ is a cone, typically a 
Cartesian product of simple, standard cones,
such as the nonnegative cone, second-order cone, and exponential cone.
Such a representation allows $f$ to be used in any DSL based
on DCP, in the objective or constraint functions.

Specifically we work out a graph form for the perspective 
\[
P(w,\delta) = \delta K(w/\delta),
\]
where $K=S(g_1, \ldots, g_k)$, with 
\[
g_i(w) = \log(\pi_i)-\mu_i^Tw+\frac{1}{2}w^T\Sigma_i w, \quad i=1, \ldots, k,
\]
and $S$ is the soft-max or log-sum-exp function \eqref{e-log-sum-exp}.

\subsection{Graph form calculus}

We view $P$ as a composition of four operations:
an affine pre-composition, then an affine post-composition, then composition, and
finally, the perspective.
We show here generic methods for carrying out these operations using 
graph form representations.
The first three operations, affine pre-composition, affine post-composition,
and composition, are known (and indeed, used in all DSLs for convex optimization);
we give them here for completeness.
The last one, the perspective transform, is not well known, but is mentioned
in \cite{persp}.

\paragraph{Affine pre-composition.}
Suppose $f$ has graph form
\[
\epi f = \{ (x,t) \mid\exists z~ Fx+Gz+td+e \in C\},
\]
and $g$ is the affine pre-composition $g(x)=f(Ax+b)$.
Then $g$ has graph form
\begin{equation}\label{e-affine-pre}
\epi g =\{(x,t)\mid \exists z~FAx+Gz+td+(Fb+e)\in C\}.
\end{equation}

\paragraph{Affine post-composition.}
Suppose $f$ has graph form
\[
\epi f = \{ (x,t) \mid\exists z~ Fx+Gz+td+e \in C\},
\]
and $h(x)=af(x)+b$, where $a\in\reals_+$ and $b\in\reals$.
Then $h$ has graph form
\begin{equation}\label{e-affine-map}
\epi h = \{(x,t)\mid\exists z~Fx+Gz+t(d/a)+(e-(b/a)d)\}.
\end{equation}

\paragraph{Composition.}
Suppose that $g_i$ are convex functions with graph forms
\[
\epi g_i = \{ (w,t_i) \mid\exists z_i~ F_iw+G_iz_i+t_id_i+e_i \in C_i\},
\quad i=1, \ldots, k,
\]
and $S$ is a convex function with graph form
\[
\epi S = \{ (u,t) \mid\exists z_0~ F_0u+G_0z_0+td_0+e_0 \in C_0\}.
\]
We assume that $S$ is increasing in each of its arguments, so the 
composition $K=S(g_1, \ldots, g_k)$ is convex.
Then $K$ has graph form
\begin{equation}\label{e-comp}
\epi K = \left\{(w, t) ~\middle |
\begin{array}{l}
\exists z_0,t_i,z_i~F_0(t_1,\ldots,t_k)+G_0z_0+td_0+e_0 \in C_0,\\
F_iw+Gz_i+t_id_i+e_i \in C_i, ~i=1, \ldots, k
\end{array}\right\}.
\end{equation}
(We can stack the affine functions of $(w,t)$, and use the product cone $C_0 \times \cdots 
\times C_k$ as the cone in the representation of $K$.)

\paragraph{Perspective.}
Here we show how to construct a graph form of the perspective
of a function given in graph form.
The perspective of $f:\reals^n\to\reals\cup\{\infty\}$ is the function
$p:\reals^{n+1}\to\reals\cup\{\infty\}$ defined by
\[
p(x,s) = \left\{
\begin{array}{ll}
  sf(x/s) & s>0\\
  0 & s=0,~x=0\\
  \infty & \mbox{otherwise.}
\end{array}\right.
\]
(See \cite[\S 3.2.6]{cvxbook}) or \cite[\S IV.2.2]{urruty1993convex}.)
Then $p$ has graph form given by 
\[
\epi p = \{(x,s,t) \mid sf(x/s)\leq t\}=\{(x,s,t)\mid f(x/s)\leq t/s\}.
\]
Substituting this expression into the graph form of $f$ given in
\eqref{e-general-epigraph}, we have
\[
\epi p = \{(x,s,t) \mid\exists z~ F(x/s)+Gz+d(t/s)+e \in C\}.
\]
Since $s>0$ and $C$ is a cone, we have
\BEAS
F(x/s)+Gz+d(t/s)+e \in C &\iff& Fx+G(sz)+dt+se \in C.
\EEAS
Thus, introducing a new affine description
\[
\tilde F = \left[\begin{array}{ccc}
F & e 
\end{array}\right],\quad \tilde e = 0,
\]
and with $\tilde{z}$ a new auxilliary variable,
$p$ has graph form
\begin{equation}\label{e-perspective}
\epi p = \{(x,s,t)\mid\exists\tilde{z}~\tilde{F}(x,s)+G\tilde z+td +\tilde e \in C\}.
\end{equation}

\paragraph{Graph form of log-sum-exp.}
\begin{eqnarray*}
S(x) \leq t &\iff& 
\log\left(\sum_{i=1}^k\exp(x_i)\right) \leq t \\
&\iff& \log\left(\sum_{i=1}^k\exp(x_i-t)\right)\leq 0\\
&\iff& \sum_{i=1}^k\exp(x_i-t)\leq 1\\
&\iff& \sum_{i=1}^ku_i\leq 1,\qquad (x_i-t,1,u_i)\in C_{\exp},\quad i=1,\ldots,k,
\end{eqnarray*}
where 
\[
C_{\exp}= \{(a,b,c)\mid e^{a/b}\leq c/b,~b> 0\}\cup\{(a,0,c)\mid a\leq 0,~c\geq 0\}
\]
is the exponential cone \cite{glineur2000,chandrasekaran2017}, which is supported
by several solvers. 
So $S$ has graph form given by
\begin{equation}\label{e-lse-cone}
\epi S = \{(x,t)\mid\exists u~ F^\text{LSE}x+G^\text{LSE}u+td^\text{LSE}+e^\text{LSE}\in
C^\text{LSE}\},
\end{equation}
with
\[
F^\text{LSE}=\left[\begin{array}{c}
0      \\ \hline
e_1^T    \\
0      \\
0      \\ \hline
\vdots \\ \hline
e_k^T    \\
0      \\
0      
\end{array}\right],\qquad
G^\text{LSE}=\left[\begin{array}{c}
\ones^T\\ \hline
0\\
0\\
e_1^T\\ \hline
\vdots\\ \hline
0\\
0\\
e_k^T\\ 
\end{array}\right],\qquad
d^\text{LSE}=\left[\begin{array}{r}
0\\ \hline
-1\\
0\\
0\\ \hline
\vdots\\ \hline
-1\\
0\\
0\\
\end{array}\right],\qquad 
e^\text{LSE}=\left[\begin{array}{r}
-1\\ \hline
0\\
1\\
0\\ \hline
\vdots\\ \hline
0\\
1\\
0\\
\end{array}\right],
\]
\[
F^\text{LSE}\in \reals^{(3k+1)\times k},\qquad G^\text{LSE}\in\reals^{(3k+1)\times k},\qquad 
d^\text{LSE}\in\reals^{3k+1},\qquad e^\text{LSE}\in\reals^{3k+1},
\]
and 
\[
C^\text{LSE} = \reals_-\times C_{\exp}\times\cdots \times C_{\exp}.
\]
The horizontal dividers denote separate blocks. After the first row, blocks of
size 3 are repeated $k$ times.
\paragraph{Graph form of quadratic.}
To derive a graph form for the function $f(x)=x^Tx$ with $x\in\reals^n$, we first observe that
\[x^Tx\leq t\iff 
\left\|
\left[\begin{array}{c}
x\\
\frac{t-1}{2}
\end{array}\right]
\right\|_2
\leq \frac{t+1}{2}.
\]
Therefore,
\begin{equation}\label{e-quad-cone}
\epi f = \{(x,t)\mid F^\text{quad}x+
td^\text{quad}+
e^\text{quad}\in C_\text{SOCP}\},
\end{equation}
with $F^\text{quad}\in\reals^{(n+2)\times n},~d^\text{quad}\in\reals^{n+2},~e^\text{quad}\in\reals^{n+2}$
defined by
\[
F^\text{quad}=\left[\begin{array}{c}
I\\
0\\
0\\
\end{array}\right],\qquad
d^\text{quad}=\left[\begin{array}{c}
0\\
1/2\\
1/2\\
\end{array}\right],\qquad
e^\text{quad} = \left[\begin{array}{c}
0\\
-1/2\\
1/2\\
\end{array}\right],
\]
and where $C_\text{SOCP}=\{(x,t)\mid \|x\|_2\leq t\}$ is the second order cone 
\cite[\S4.4.2]{cvxbook}\cite{nesterov1994interior}.

\subsection{Graph form of EVaR}
Using the calculus outlined above, we can now develop a graph form
of $P$, where $P(w,\delta) = \delta K(w/\delta)$.
First, we use affine pre-composition to write 
\[
g_i(w) = \log(\pi_i)-\mu_i^Tw+\frac{1}{2}w^T\Sigma_i w\leq t
\]
as
\[
f(A_iw+b_i)-\frac{1}{2}\mu_i^T\Sigma_i^{-1}\mu_i+\log(\pi_i)\leq t,
\] 
where
\[
f(w)=w^Tw,\qquad A_i=\frac{1}{\sqrt{2}}\Sigma_i^{1/2},\qquad 
b_i=-\frac{\sqrt{2}}{2}\Sigma_i^{-1/2}\mu_i.
\] 
Thus, using our affine pre-composition expression \eqref{e-affine-pre} together
with our graph form of the quadratic \eqref{e-quad-cone} and affine 
post-composition \eqref{e-affine-map}, we have
\[
\epi g_i = \left\{(w,t_i)~\middle|~ \exists z_i~
(F^\text{quad}A_i)w+
t_id^\text{quad}+
e_i^\text{quad}\in C_\text{SOCP}\right\},
\]
with
\[
e_i^\text{quad} = F^\text{quad}b_i+e^\text{quad}+
\left(\frac{1}{2}\mu_i^T\Sigma_i^{-1}\mu_i-\log(\pi_i)\right)d^\text{quad}.
\]
Then, using the graph form of log-sum-exp given in \eqref{e-lse-cone} and the
composition rule give in \eqref{e-comp}, we can write the composition in graph
form as
\BEAS
\epi K &=& \left\{(w, t) ~\middle |~\exists z_0,t_1,\ldots,t_k\quad
\begin{array}{l}
F^\text{LSE}(t_1,\ldots,t_k)+G^\text{LSE}z_0+td^\text{LSE}+e^\text{LSE} \in C^\text{LSE},\\
F^\text{quad}A_iw+t_id^\text{quad}+ e_i^\text{quad}\in C_\text{SOCP}, ~i=1, \ldots, k
\end{array}\right\}\\
&=& \{(w,t)\mid \exists z~F^Kw+G^Kz+td^K+e^K\in C^K\},
\EEAS
with
$F^K\in\reals^{(3k+1+k(k+2))\times k},~G^K\in\reals^{(3k+1+k(n+2))\times 2k},
~d^K\in\reals^{3k+1+k(n+2)},~e^K\in\reals^{3k+1+k(n+2)}$
defined by
\[
F^K = \left[\begin{array}{c}
0\\
F^\text{quad}A_1\\
\vdots\\
F^\text{quad}A_k
\end{array}\right],\qquad 
G^K = \left[\begin{array}{cc}
G^\text{LSE} & F^\text{LSE}\\
0 & d^\text{quad}e_1^T\\
\vdots & \vdots \\
0 & d^\text{quad}e_k^T \\
\end{array}\right],\qquad 
d^K = \left[\begin{array}{c}
d^\text{LSE}\\0\\\vdots\\0\\
\end{array}\right],\qquad 
e^K = \left[\begin{array}{c}
e^\text{LSE}\\e_1^\text{quad}\\\vdots\\e_k^\text{quad}
\end{array}\right],
\]
and
\[
C^K = C^\text{LSE}\times C^\text{SOCP}\times \cdots \times C^\text{SOCP}.
\]
Finally, using the perspective rule given in \eqref{e-perspective},
the perspective of $K$ has graph form given by
\[
\{(w,\delta,t)\mid \delta K(w/\delta)\leq t\} =
\{(w,\delta,t)\mid \exists z~\tilde F^K(w,\delta)+G^Kz+td^K\in C^K\},
\]
with 
$\tilde F^K = \left[\begin{array}{cc}
F^K & e^K
\end{array}\right]$.

\paragraph{CVXPY specification.}
CVXPY code for EVaR portfolio optimization using the graph form of $\delta K(w/\delta)$
is available at the repository
\url{https://github.com/cvxgrp/exp_util_gm_portfolio_opt}.

\end{document}